\documentclass[conference]{IEEEtran}
\IEEEoverridecommandlockouts

\usepackage{cite}
\usepackage{amsmath,amssymb,amsfonts}
\usepackage{graphicx}
\usepackage{textcomp}
\usepackage{xcolor}
\usepackage{xspace}
\usepackage{booktabs}
\usepackage{algorithm}
\usepackage{algpseudocode}
\usepackage{xurl}
\usepackage{tabularx}
\usepackage{siunitx}
\usepackage{comment}
\allowdisplaybreaks[1]

\definecolor{lightgray}{gray}{0.9}
\definecolor{darkergray}{gray}{0.8}
\definecolor{forestgreen}{RGB}{34,139,34}

\usepackage{hyperref}
\usepackage[capitalize]{cleveref}
\hypersetup{
    colorlinks=true,
    citecolor=green,
    filecolor=black,
    linkcolor=red,
    urlcolor=blue
}


\newcommand{\B}{\mathcal{B}}

\newcommand{\Z}{\mathcal{Z}}
\newcommand{\N}{\mathcal{N}}
\newcommand{\G}{\mathcal{G}}

\title{\LARGE \bf 
A Fast and Accurate Time-Space Decomposition\\ for Unit Commitment With Storage
}

\author{Chance Brainard$^{1,\dagger}$, Kylind Reagan$^{1,\dagger}$, Shaked Regev$^{1,\dagger}$% <-this % stops a space
\thanks{Notice: This manuscript has been authored by UT-Battelle, LLC, under contract DE-AC05-00OR22725 with the US Department of Energy (DOE). The US government retains and the publisher, by accepting the article for publication, acknowledges that the US government retains a nonexclusive, paid-up, irrevocable, worldwide license to publish or reproduce the published form of this manuscript, or allow others to do so, for US government purposes. DOE will provide public access to these results of federally sponsored research in accordance with the DOE Public Access Plan (\href{https://www.energy.gov/doe-public-access-plan}{https://www.energy.gov/doe-public-access-plan}). Research sponsored by the Laboratory Directed Research and Development Program of Oak Ridge National Laboratory, managed by UT-Battelle, LLC, for the US Department of Energy.}% <-this % stops a space
\thanks{$^{1}$Oak Ridge National Laboratory, 1 Bethel Valley Road, Oak Ridge, Tennessee, USA.
        {\tt\small regevs@ornl.gov}}%
\thanks{$^{\dagger}$Equal Contribution}%
}
\date{October 2026}

\begin{document}

\maketitle

\begin{abstract}
Grid-scale batteries couple decisions across time. This makes unit commitment (UC) with thousands of  buses, batteries, generators, and alternating current constraints intractable.
We introduce a time-space decomposition (TSD) heuristic to solve the problem. We first calculate target net battery charge or discharge rates on a zonal level, relaxing spatial constraints. We then incorporate this target into the Relax-and-Round UC (RRUC) solver from our previous work.  RRUC enforces the spatial transmission constraints and transfers power between zones, if that reduces costs. We show that our TSD produces an error of about 0.1\% on a simpler economic dispatch problem with batteries. This is larger than the error produced by RRUC in the tests in our previous work. So, we conjecture that the combined error of our two approximations is not much larger than 0.1\%. We show that TSD with RRUC can solve systems with 10,000 buses in about 6 minutes, with sub-quadratic scaling, making it suitable for real time operations on an ISO-sized system.

% develop a solver present a solver with two capabilities: optimizing batteries over multiple time intervals and committing generators with Relax-and-Round UC (RRUC), which previously produced near-optimal commitments in near linear time.
% A rolling-horizon quadratic program (QP) plans each zone's battery target and an optimal power flow produces a zonal marginal price that turns the target into an objective term, so the AC solver is free to place batteries on the network.
% We combine these methods into one solver we call Time Space Decomposed RRUC. We show that this method scales as $\mathcal{O}(n^{1.31})$. Because each capability is validated separately, with the QP within $\mathbf{0.12\%}$ of the time-coupled solution's cost on an economic dispatch, this combination can be trusted where no time-coupled UC benchmark exists.
\end{abstract}


\section{Introduction}
\label{sec:Intro}
Energy storage is an effective way to balance the electric grid when electricity demand is particularly high or low relative to supply. It can react quickly to changing grid conditions and can operate without supporting electricity~\cite{css_us_energy_storage_factsheet_2026}. Storage
accounts for 5-6\% of the total electric grid capacity in the United States~\cite{eia_todayinenergy_67925,css_us_energy_storage_factsheet_2026} and 2-3\% worldwide~\cite{hydropower2026outlook}. 

Traditional generators can continuously generate electric power for long periods of time as long as they continue to have fuel, their external source of energy. Storage however, produces net zero electric energy. Storage operation couples every period to the next through its state of charge (\textbf{SoC}). This means that loose time coupling approaches,  which ignore intermediate states and only look at peaks (as in \cite{regev2026economicUC,brody2026fastrruc}), are not suitable for battery optimization.

Unit Commitment (\textbf{UC}) with storage is conventionally posed as a mixed-integer program (\textbf{MIP}) in which SoC dynamics link all periods of the horizon~\cite{achterberg2013}.
This adds continuous variables and linking constraints to a problem that is already NP-hard, and with transmission modeled by piecewise-linear alternating current (\textbf{AC}) approximations the formulation scales poorly~\cite{lawler1966, morrison2016}. 
This means that time-coupled UC with storage constraints cannot be solved at the scale of modern grids.

Relax-and-Round UC (\textbf{RRUC})~\cite{regev2026economicUC,regev2026oneperiodUC, brody2026fastrruc} avoids directly solving a MIP. First, it relaxes the binary commitment variables to $[0,1]$ and solves the relaxed problem with a nonlinear interior-point solver (\textbf{IPS}). Next, it sorts units by their relaxed value and constructs different sets of generators that can feasibly supply demand. Finally, it solves the corresponding AC optimal power flow (\textbf{ACOPF}) problems with transmission constraints to verify feasibility and choose the solution with the best objective.

We extend this formulation to storage. Lithium-ion batteries and hydroelectric pumped storage are the two types of energy storage. Here, we will not account for the differences between them and refer to them collectively as `batteries'. We solve the intertemporal part of the problem by a cheap, structured look-ahead, leaving the per-period ACOPF UC problem single-period. Our main contributions in this paper are:
\begin{itemize}
    \item Introducing a time-space decomposition  extension to RRUC (\textbf{TSD-RRUC}) which can commit batteries over multiple periods.
    \item Showing TSD-RRUC scales as $\mathcal{O}(n^{1.31})$ and can solve systems with 10,000 buses in about six minutes.
    \item Showing that the TSD heuristic has an error of about 0.1\% on our economic dispatch (\textbf{ED}) test cases\footnote{\cite{regev2026economicUC,brody2026fastrruc} show that RRUC's heuristic has even smaller errors. This makes it plausible that the combined heuristic has an error not much larger than 0.1\%.}. 
\end{itemize}

\cref{sec:formulation} details our problem formulation, solution method, and the steps of our algorithm. \cref{sec:design} describes generators for scalable battery fleets and synthetic zone demand. \cref{sec:results} compares cost and runtime against the time-coupled solution for ED and mean-per-period RRUC solve time. \cref{sec:conclusion} summarizes our results. 

\section{Formulation and Methods}
\label{sec:formulation}
\subsection{Mathematical Formulation}
\label{subsec:math}
Our mathematical formulation of UC extends our previous formulation to include batteries~\cite{brody2026fastrruc} (see for further details). We include an overview here. \cref{tab:parametersets} describes the input parameter sets of our UC optimization problem. 

\begin{table}[htbp]
    \centering
    \caption{Input parameter sets for the optimization solver}
    \label{tab:parametersets}
    
    \renewcommand{\arraystretch}{1.1}
    \noindent
    \begin{tabularx}{0.5\textwidth}{@{} c c X @{}}
    \toprule
    \textbf{Notation} & \textbf{Indexing} &\multicolumn{1}{c}{\textbf{Meaning}} \\\midrule
    \(\mathcal{N}\) & & Set of buses \\\midrule
    \(\mathcal{G},\mathcal{B}\subseteq\mathcal{N}\) & & Set of generators $\mathcal{G}$ and set of batteries $\mathcal{B}$\\\midrule
    \(\mathcal{G}_m, \mathcal{G}_d \subseteq\mathcal{G}\) & & Set of discretionary and must-run generators \\\midrule
    \(\mathcal{S}\) & & Set of contingency scenarios \\\midrule
    \(\mathcal{L},N(i)\) & & Set of lines and set of buses connected to bus \(i\) \\\midrule
    \(\mathcal{L}_\delta,N_\delta(i)\) & \(\delta\in\mathcal{S},\) & Set of lines and set of buses connected to bus \(i\) in contingency \(\delta\) \\\bottomrule
    \end{tabularx}
\end{table}

\cref{tab:decision variables} describes the decision variables of our problem. The $u_i$s, $P_i$s, and $\rho_i$s are the variables the decision maker would need to commit units. The other variables are the result of the underlying physics. Variables with subscript \(\delta\) are the values of the corresponding variables during contingency \(\delta\). 

\begin{table}[htbp]
    \centering
    \caption{Decision Variables for the optimization solver}
    \label{tab:decision variables}
    
    \renewcommand{\arraystretch}{1.1}
    \noindent
    \begin{tabularx}{0.5\textwidth}{@{} c c X @{}}
    \toprule
    \textbf{Notation} & \textbf{Indexing} &\multicolumn{1}{c}{\textbf{Meaning}} \\\midrule
    \(u_i\) & \(i\in\mathcal{G}_d\) & Commitment status \\\midrule
    \(P_i,Q_i\) & \(i\in\mathcal{G}\) & Real and reactive power production \\\midrule
    \(\rho_i, C_i\) & \(i\in\mathcal{B}\) & Battery power flow and capacity \\\midrule
    \(V_i,\theta_i\) & \(i\in\mathcal{N}\) & Voltage magnitude and phase angle 
    \\\midrule
    \(P_{ij},Q_{ij}\) & \(ij\in\mathcal{L}\) & Real and reactive power from bus \(i\) to bus \(j\) as measured at bus \(i\) \\\midrule
    \(P_{R,\delta,i},Q_{R,\delta,i}\) & \(i\in\mathcal{N}, \delta\in\mathcal{S}\) & Real and reactive power from a generator that is ramping up at bus \(i\), during contingency \(\delta\) \\\bottomrule
    \end{tabularx}
\end{table}

Along with the input parameters defined in Table 3 of~\cite{brody2026fastrruc}, we introduce additional input parameters for batteries indexed by $i\in \mathcal{B}$. \(\rho_i^{\max}\) is the battery real power flow limit and \(C_i^{\min}, C^{\max}_i\) are the battery capacity limits.

With these we can define our optimization problem:

\begin{subequations}
\label{eq:UnitCommitment}
\begin{align}
\label{eq:Objective}
\min \; &
\sum_{i \in \mathcal{G}_m} \left(a_i P_i^2 + b_i P_i + c_i\right) 
\nonumber\\&+ \sum_{i \in \mathcal{G}_d} \left(a_i P_i^2 + b_i P_i + c_i\right) u_iu_{t-1,i}
\nonumber\\&+ \sum_{i \in \mathcal{G}_d}K_i(u_{t-1,i}(1-u_i)+(1-u_{t-1,i})u_i) \nonumber\\
& + \gamma\sum_{\delta\in\mathcal{S}}\sum_{i\in \mathcal{G}}u_i\left(a_i P_{\delta,i}^2 + b_i P_{\delta,i} + c_i\right)
\\[6pt]
\text{s.t.} \;
& \forall i\in\mathcal{G} : \max\!\left(P_{\min,i},\, P_{t-1,i} - r_{d,i}\right) \le P_i \nonumber\\&
\le \min\!\left(P_{\max,i},\, P_{t-1,i} + r_{u,i}\right)
\label{eq:ConsRamp}
\\& 
\forall i \in \B : 
C_i^{\min} \le C_i - \rho_i \le C_i^{\max} 
\label{eq:chargeBatt}
\\& 
\forall i \in \B : 
-\rho_i^{\max} \le \rho_i \le \rho_i^{\max}
\label{eq:rateBatt}
\\&
\forall i\in\mathcal{G} : Q_{\min,i} \leq Q_i \leq Q_{\max, i}
\label{eq:ConsMVar}
\\&
\forall i\in\mathcal{G} : u_i\in \{0,1\}
\label{eq:ConsBools}
\\&
\forall i \in \mathcal{N} : \;(P_{i}u_i + P_{\min,i}(1-u_i))u_{t-1,i} +\nonumber\\& P_{R,i} + \rho_i - P_{Di}  =G_{\mathrm{sh},i} + \sum_{j \in N(i)} P_{ij} \nonumber\\& \doteq
\sum_{j \in N(i)\cup\{i\}}  V_{i}V_{j}
\left( G_{ij}\cos\theta_{ij} + B_{ij}\sin\theta_{ij} \right)
\label{eq:ConsPflow}
\\&
\forall i \in \mathcal{N} : \;Q_{i}u_iu_{t-1,i} - Q_{Di}=-B_{\text{sh},i} + \sum_{j \in N(i)} Q_{ij} \nonumber\\& \doteq 
\sum_{j \in N(i)\cup\{i\}}V_{i}  V_{j}
\left( G_{ij}\sin\theta_{ij} - B_{ij}\cos\theta_{ij} \right)
\label{eq:ConsQflow}
\\&
\forall ij\in\mathcal{L} : P_{ij}^2 + Q_{ij}^2 \le S_{\max,ij}^2
\label{eq:ConsThermal}
\\&
\forall i\in\mathcal{N} : V_{\min,i} \le V_{i} \le V_{\max,i}
\label{eq:ConsVoltage}
\\&
\forall ij\in\mathcal{L} : \theta_{ij} = \theta_{i} - \theta_{j}
\label{eq:ConsAngle}
\\&
\forall i\in\mathcal{G} : P_{\min,i} \leq P_{\delta,i},P_{R,\delta,i} \leq P_{\max,i}
\label{eq:ConsPDelta}
\\&
\forall i\in\mathcal{G} : Q_{\min,i} \leq Q_{\delta,i},Q_{R,\delta,i} \leq Q_{\max,i}
\label{eq:ConsQDelta}
\\& 
\forall i \in \B : 
-\rho_i^{\max} \le \rho_{\delta,i} \le \rho_i^{\max}
\label{eq:rateBattDelta}
\\&
\forall i \in \mathcal{N}, \;\delta\in \mathcal{S}: \; P_{\delta,i}u_i + P_{R,\delta,i} + \rho_{\delta,i} - P_{\delta,Di}\nonumber = \nonumber\\&G_{\text{sh},i}+\sum_{j \in N_\delta(i)} P_{\delta,ij} \nonumber\\& \doteq
\sum_{j \in N_\delta(i)\cup\{i\}}  V_{\delta,i}V_{\delta,j}
\left( G_{ij}\cos\theta_{\delta,ij} + B_{ij}\sin\theta_{\delta,ij} \right)
\label{eq:ConsPflowdel}
\\&
\forall i \in \mathcal{N}, \delta\in \mathcal{S}: \; Q_{\delta,i}u_i + Q_{R,\delta,i} - Q_{\delta,Di}\nonumber\\& = -B_{\text{sh},i}+\sum_{j \in N_\delta(i)} Q_{\delta,ij} \nonumber\\& \doteq
\sum_{j \in N_\delta(i)\cup\{i\}}V_{\delta,i}  V_{\delta,j}
\left( G_{ij}\sin\theta_{\delta,ij} - B_{ij}\cos\theta_{\delta,ij} \right)
\label{eq:ConsQflowdel}
\\&
\forall \delta\in\mathcal{S}, ij\in\mathcal{L}_\delta : P_{\delta,ij}^2 + Q_{\delta,ij}^2 \le S_{\max,ij}^2
\label{eq:ConsThermaldel}
\\&
\forall \delta\in\mathcal{S}, i\in\mathcal{N} : V_{\min,i} \le V_{\delta,i} \le V_{\max,i}
\label{eq:ConsVoltagedel}
\\&
\forall \delta\in\mathcal{S}, ij\in\mathcal{L}_\delta : \theta_{\delta,ij} = \theta_{\delta,i} - \theta_{\delta,j}
\label{eq:ConsAngledel}
\end{align}
\end{subequations}

Storage adds the battery constraints~\cref{eq:chargeBatt,eq:rateBatt,eq:rateBattDelta} to the storage-free formulation in~\cite{brody2026fastrruc}. It also adds the $+ \rho_i$ in~\cref{eq:ConsPflow} and the $+ \rho_{\delta,i}$ in~\cref{eq:ConsPflowdel}.

This problem becomes even more complicated because we must optimize \cref{eq:UnitCommitment} over multiple time periods simultaneously, which is intractable. Instead, we introduce a time-space decomposition. We will first calculate a target zone-level battery output that reduces the variance of the demand. This reduces total productions costs because it grows super-linearly with demand~\cite{regev2026fastaccuratehydrounit}. We then solve a single period UC with an added term that `encourages' batteries to charge or discharge based on the multi-period plan from the previous step.

\subsection{Zonal Targets}
\label{subsec:zone-targ}
Our system is subdivided into \textit{zones}, while maintaining MVA constraints on each transmission element so that the model remains nodal. Each resulting zone has relatively high internal connectivity and weaker connectivity to the rest of the system. See \cite{brody2026fastrruc,tsybina2026agenticartificialintelligencepower} for more details about our congestible and extendable system design.

We reach a tentative battery schedule by solving a problem to reduce demand variance in time within each zone. At period $t_0$, given a demand forecast $d_{i,t}$ for each bus $i$ over the window $[t_0, T]$, our formulation chooses battery flows $\rho_{i,t}$ that flatten the zone's net demand, subject to the batteries' power, energy, and terminal-SoC requirement (which we set to the initial SoC, without loss of generality).

\begin{equation}
\label{eq:target}
\begin{aligned}
\min \quad & \sum_{t = t_0}^{T} \Bigl( \sum_{i\in \N_z} d_{i,t} - \sum_{i \in \B_z} \rho_{i,t} \Bigr)^{2} \\[6pt]
\text{s.t.} \;
& \forall\, i \in \B _z, \; \forall t \in t_0,...,T:\\
&  C_{i,t} = C_{i,t-1} - \rho_{i,t}; \quad C_{i,T} = C_{i,t_0} \\
&  -\rho_i^{\max} \le \rho_{i,t} \le \rho_i^{\max}; \quad C_i^{\min} \le C_{i,t} \le C_i^{\max} 
\end{aligned}
\end{equation}
where $\N_z$ and $\B_z$ are the buses and batteries of zone $z$. We use the $\rho_i^*$s to construct $R^* = \sum_{i\in\mathcal{B}} \rho_i^*$ as the optimal net battery flow for current period (before constraints).

\subsection{Marginal Cost}
\label{subsec:mc}
The marginal cost of a system is its cost to produce slightly more or less power.
We compute marginal cost by solving the economic dispatch (\textbf{ED}) problem on the current period using must-run and discretionary on generators only. Mathematically, we set $u_i=1$ for  on and $u_i=0$ for off discretionary generators. We do so while allowing batteries in each zone to have summed flow $\rho_i$ equal to that of $R^*$ for that zone as computed in \cref{subsec:zone-targ}. Using the solution, we compute the dual of generation cost with respect to net battery power flow for each zone. This gives us the marginal cost $MC^*_z$ for each zone when battery power is discharged or charged according to variance reduction optimum.

\subsection{Algorithm}
\label{subsec:algo}
In the final step of our algorithm, we incorporate information about future time periods. We do this by adding the term $\sum_{z\in \Z} MC^*_z\sum_{i\in\B_z}\rho_i$ to \cref{eq:Objective}, where $\Z$ is the set of zones. This keeps the intertemporal coupling within one term, making solving the modified \cref{eq:UnitCommitment} tractable. 
% \begin{equation}
%     \label{eq:MC-term}
    
% \end{equation}

\cref{alg:storage-rruc} details a single period of the solver. It builds on the algorithm of~\cite{brody2026fastrruc}, which we apply to slightly different problem (Steps 6 and 7 in \cref{alg:storage-rruc}).

\begin{algorithm}[H]
\caption{One period of TSD-RRUC}
\label{alg:storage-rruc}
\begin{algorithmic}[1]
    \Require UC state, battery SoC $C_i$, forecast window $[t_0, T]$
    \For{each zone $z\in\Z$}
        \State $\rho^*_{i,t}\gets$ solution of~\cref{eq:target} for $i\in\B_z$
        \State $R^*_z\gets\sum_{i\in\B_z}\rho^*_{i, t_0}$.
        \State $MC^*_z\gets$ dual of generation cost with respect to net battery power flow in zone $z$ (\cref{subsec:mc}). 
    \EndFor
    \State Solve \cref{eq:UnitCommitment} with $\sum_{z\in \Z} MC^*_z\sum_{i\in\B_z}\rho_i$  added to~\cref{eq:Objective}, $u_i\in\{0,1\}\to y_i\in [0,1]$ for $y^*_i,\,\rho_i$.
    \State \textsc{RoundSolution}$(\{y^*_i\})$ as in \cite{brody2026fastrruc} with $\rho_i$ held fixed.
    \State Update generator states, ramping power, and $\G_d, \G_m$.
    \State $C_i\gets C_i-\rho_i$, $\forall i\in\B$; advance forecast  one period.
\end{algorithmic}
\end{algorithm}

\section{Test Systems}
\label{sec:design}

\subsection{Scalable Battery Cases}
\label{subsec:scale-batt}

We use the congestible 100-bus system of~\cite{brody2026fastrruc,tsybina2026agenticartificialintelligencepower} adapted from the IEEE 118-bus case~\cite{ieee118bus} to build a base case. We build a synthetic battery fleet for the 100-bus system.
% We place batteries with a heuristic that prioritizes realism over optimality, based on ERCOT~\cite{Modo-Q3-2026,ERCOT-Peak-2026,EIA-860M}. 
The fleet's total power rating is 20.5$\%$ of the system's peak demand and each battery stores 1.8 hours of energy at its rated power. 

For each zone, we place that zone's batteries at its largest-load buses, excluding buses with hydroelectric units, since hydro already behaves like storage and is rarely paired with batteries in practice~\cite{gorman2026hybrid}. We site at large loads to follow ERCOT, where battery capacity has concentrated near load centers~\cite{miller2025ercot}. Storage at these sites can serve local demand and relieve congestion thus providing most value ~\cite{thrampoulidis2013storage}. 

\begin{figure}[htbp]
    \centering
    \includegraphics[width=0.85\linewidth]{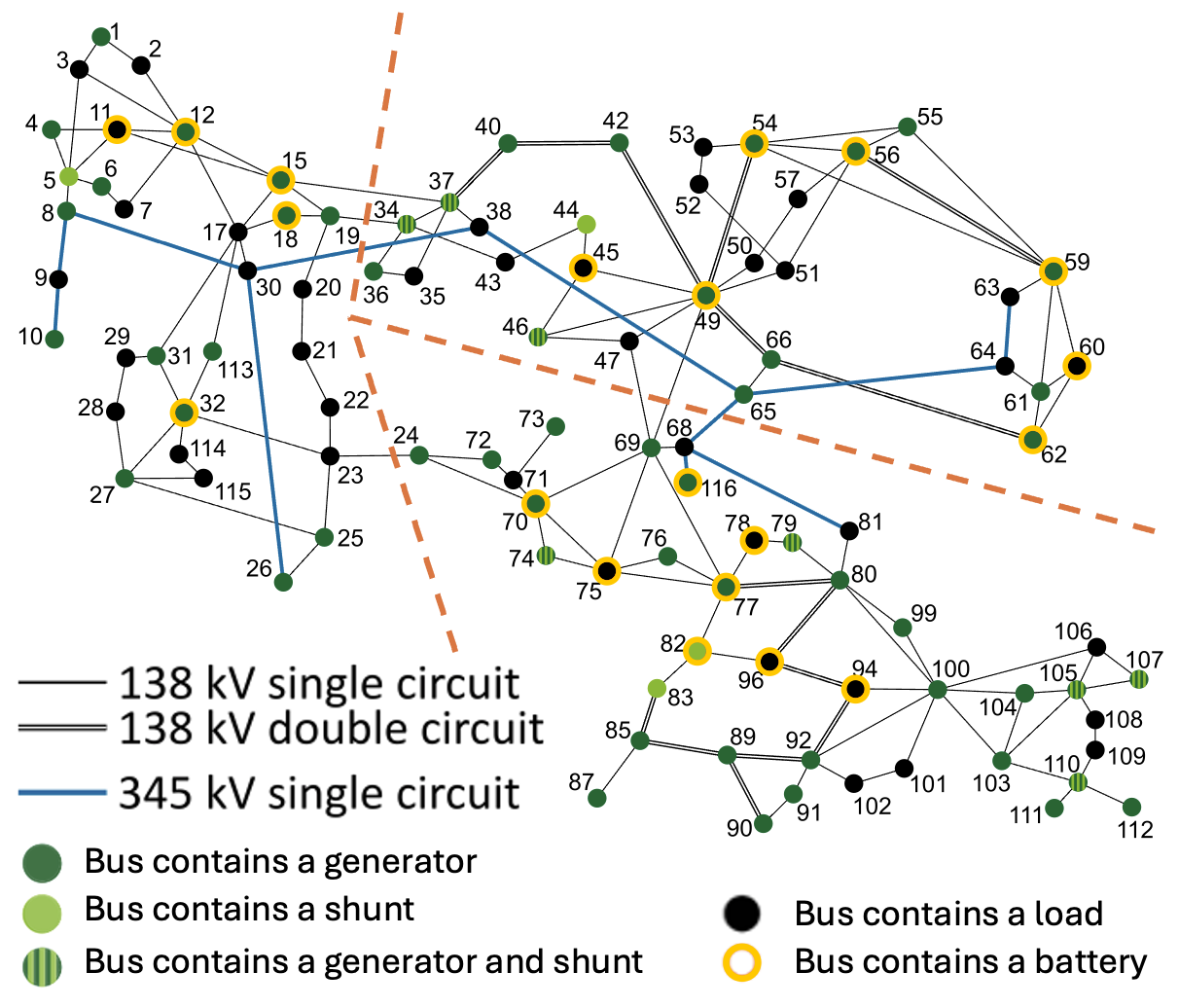}
    \caption{Modified 100 bus case based on IEEE 118 bus case~\cite{ieee118bus}. Dashed lines indicate zonal split.}
    \label{fig:grid}
\end{figure}

We place 20 batteries on the grid. We split  the battery count and the fleet's total power rating across zones in proportion to each zone's share of base case load, rounding counts to whole batteries. For the three 100-bus case zones, this gives 5, 7, and 8 batteries with 260, 394, and 479~MW of total power rating.  We divide the zone's total allocated battery power rating among these buses in proportion to their loads. Sites where the battery and any generators at its bus cannot export at rated power without a line overload, including after any single-line outage~\cite{NERC-TPL-001} are replaced by next largest load-eligible bus in their zone. ~\cref{fig:grid} illustrates the adjusted grid configuration.

%~\cite{thrampoulidis2013storage,miller2025ercot}
% 1: shows trend of ercot placements to load follow
% 2: mathematical proof for fiscal optimality of load following placements in simple system, serving as base that other papers have expanded off of, this fiscal reasoning is then the incentive for battery placements in ERCOT 

%~\cite{gorman2026hybrid}
% Batt + Hydro collocation is very rare acccording to this source
% Other sources I found worked out why it is fiscally suboptimal and nearly always avoided. The batt-hydro collocation i did find was not for traditional fiscal reasoning. Also batt-hydro collocation doesnt make much sense as hydro functions somewhat like a battery.

We produce larger bus systems by connecting perturbed duplicates of the 100-bus system randomly vertically or horizontally as in \cite{brody2026fastrruc,tsybina2026agenticartificialintelligencepower}. We perturb generator characteristics randomly up to $\pm 1\%$ and battery power flow and capacity limits up to  $\pm 2.5\%$. We place 20 batteries on each duplicate. 
%separately, repeating the full placement procedure on it as a standalone network, ignoring lines between duplicates. 
\subsection{Synthetic Zonal Demand Generation}
\label{subsec:zonal-demand}

Following~\cite{zheng2023}, we allocate zonal load to buses by fixed within-zone shares and map system load to zones via random logistic-normal weights. This allows synthetic zonal demand generation without historical zonal data.
Then each zone follows the system curve shifted by a random lag of up to 60 min, so zonal peaks no longer coincide and batteries in different zones must respond at different times.

\section{Results}
\label{sec:results}
We previously showed that on weakly time-coupled problems without batteries or transmission constraints, RRUC produces solutions that are within tolerance of Gurobi~\cite{gurobi} solutions~\cite{regev2026economicUC}. For the transmission constrained AC UC problems which Gurobi cannot solve, RRUC produces solutions within lower bounds on the optimum from relaxing the binary constraints~\cite{regev2026economicUC,brody2026fastrruc}. This supported the accuracy of our first approximation: relaxing and rounding.

An IPS is able to solve fully time-coupled ED with batteries. We will begin by showing that our space-time decomposition method produces solutions that are close to those of an IPS, on this problem. This will support the accuracy of the second of our approximations: the space-time decomposition.

We are unaware of any solver able to solve fully time-coupled UC with batteries at scale. That is to solve the problem with neither of the  approximations. This prevents us from verifying that the combination of these approximations is still a good approximation. However, it is reasonable to assume that the error this combination creates is at most as large as the sum of the errors that each approximation creates separately. With this in mind, we will present the results of our fully time-coupled AC UC battery problem without comparison.

\subsection{Economic Dispatch}
\label{subsec:results-ed}

We solve of systems of 100-2,000 buses using IPS and TSD-ED over a four-day horizon at 15 minute intervals (384 periods) across three demand scaling seeds.  

\begin{figure}[htbp]
    \centering
    \includegraphics[width=0.85\linewidth]{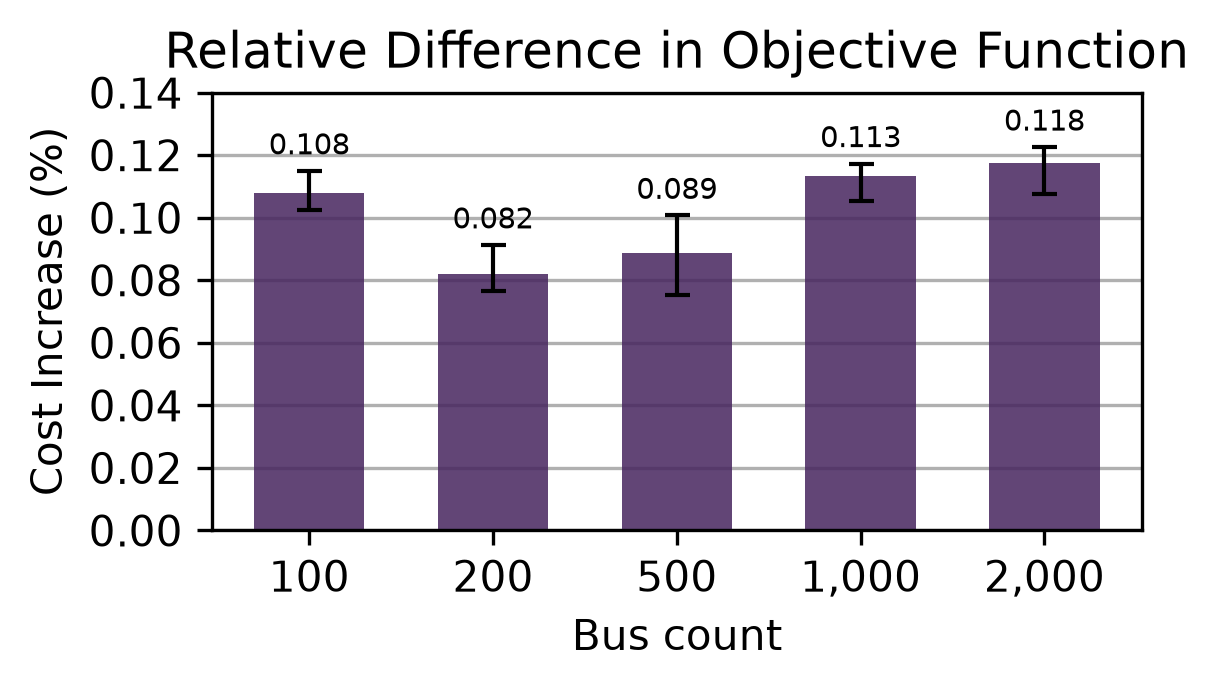}
    \caption{Experimental results for the mean cost of our TSD-ED solver relative to the fully time-coupled IPS solver over a four-day horizon (384 15-min periods) on the ED case. We solved this over 3 different demand shifts and plotted their range and mean. TSD-ED stays within about 0.12\% of the time-coupled cost at every system size.}
    \label{fig:ED_accuracy}
\end{figure}

\cref{fig:ED_accuracy} shows the relative increase in production cost incurred by using TSD-ED compared to time-coupled IPS. TSD-ED's total production cost is consistently within about 0.12\% of the time-coupled IPS solution for the systems tested. The deviations produced by different random seeds are relatively small. This supports the accuracy of the TSD.

\begin{figure}[htbp]
    \centering
    \includegraphics[width=0.85\linewidth]{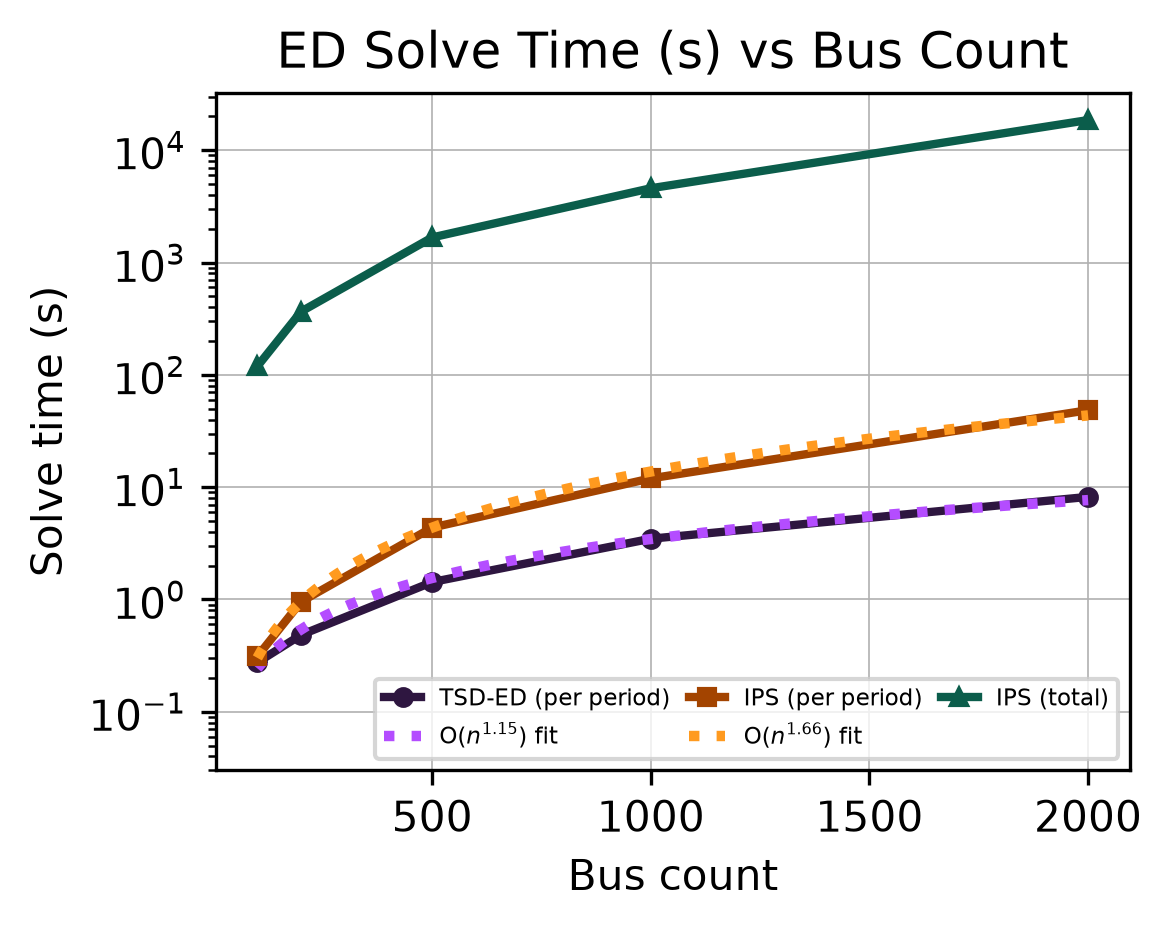}
    \caption{Results for per-period solve time versus bus count of the ED case over a four-day (384 15-min periods) horizon. Time-coupled time is its total solve time divided by total periods solved. Our TSD-ED solver scales near-linearly ($\mathcal{O}(n^{1.15})$) versus $\mathcal{O}(n^{1.66})$ for the time-coupled IPS solve.}
    \label{fig:ED_timeComplexity}
\end{figure}

\cref{fig:ED_timeComplexity} compares the run-time TSD-ED, which is naturally per period, of a period and total time-coupled IPS solve time\footnote{The per-period solve time is a good measure with perfect predictive information about the system state. One can simply solve the problem over non-overlapping time windows. The total time is a better measure if the system inputs change and the battery schedule must be reoptimized at every period.}. 
TSD-ED scales as $\mathcal{O}(n^{1.15})$, compared to time-coupled IPS's $\mathcal{O}(n^{1.66})$. TSD-ED's per-period solve is about five times faster.  

For a grid operator to make a battery dispatch decision they would need to run a fully time coupled economic dispatch solve. For 400 batteries and 2,000 buses optimized over four days, this would take about six hours. This makes optimizing systems of comparable size or larger infeasible in real time. TSD-ED achieves a similar solution in less than 10 seconds (but must be rerun for each time-period).

\subsection{Unit Commitment}
\label{subsec:results-uc}

We solve systems of 100-10,000 buses using TSD-RRUC over a four-day horizon at 15 minute intervals (384 periods).  

\begin{figure}[htbp]
    \centering
    \includegraphics[width=0.85\linewidth]{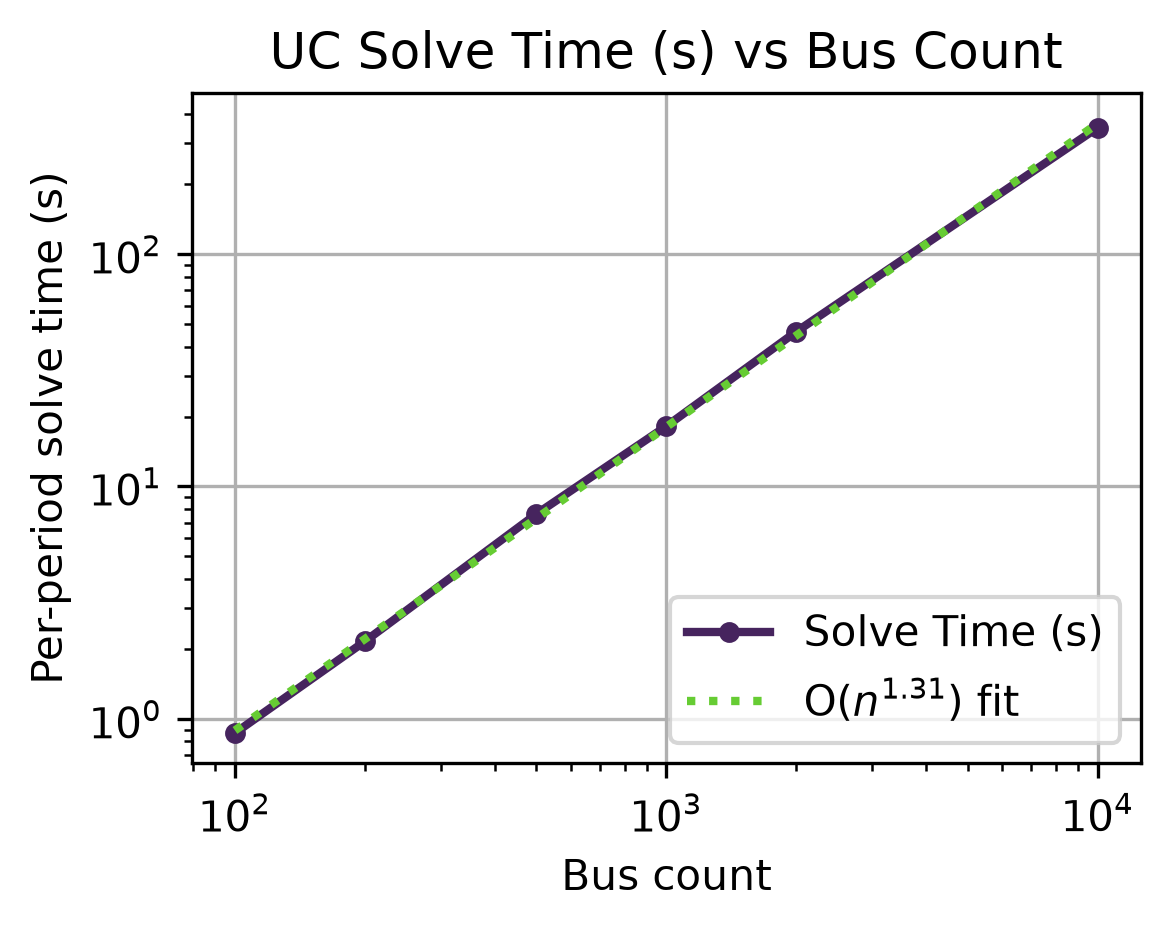}
    \caption{\cref{alg:storage-rruc}'s mean per-period solve time over a four-day (384 15-min periods) horizon. TSD-RRUC mean per-period solve grows as $\mathcal{O}(n^{1.31})$ with the number of buses.}
    \label{fig:UC_time_scaling}
\end{figure}

\cref{fig:UC_time_scaling} shows that mean per-period solve time scales as $\mathcal{O}(n^{1.31})$. A 10,000 bus system takes about 6 minutes to solve. This makes it practical to optimize systems of this size in real time, with rolling demand forecasts and grid conditions. This is similar to the scaling before adding storage of the algorithm in \cite{brody2026fastrruc}, though the overall solve time is about twice as high. 

\begin{figure}[htbp]
    \centering
    \includegraphics[width=0.85\linewidth]{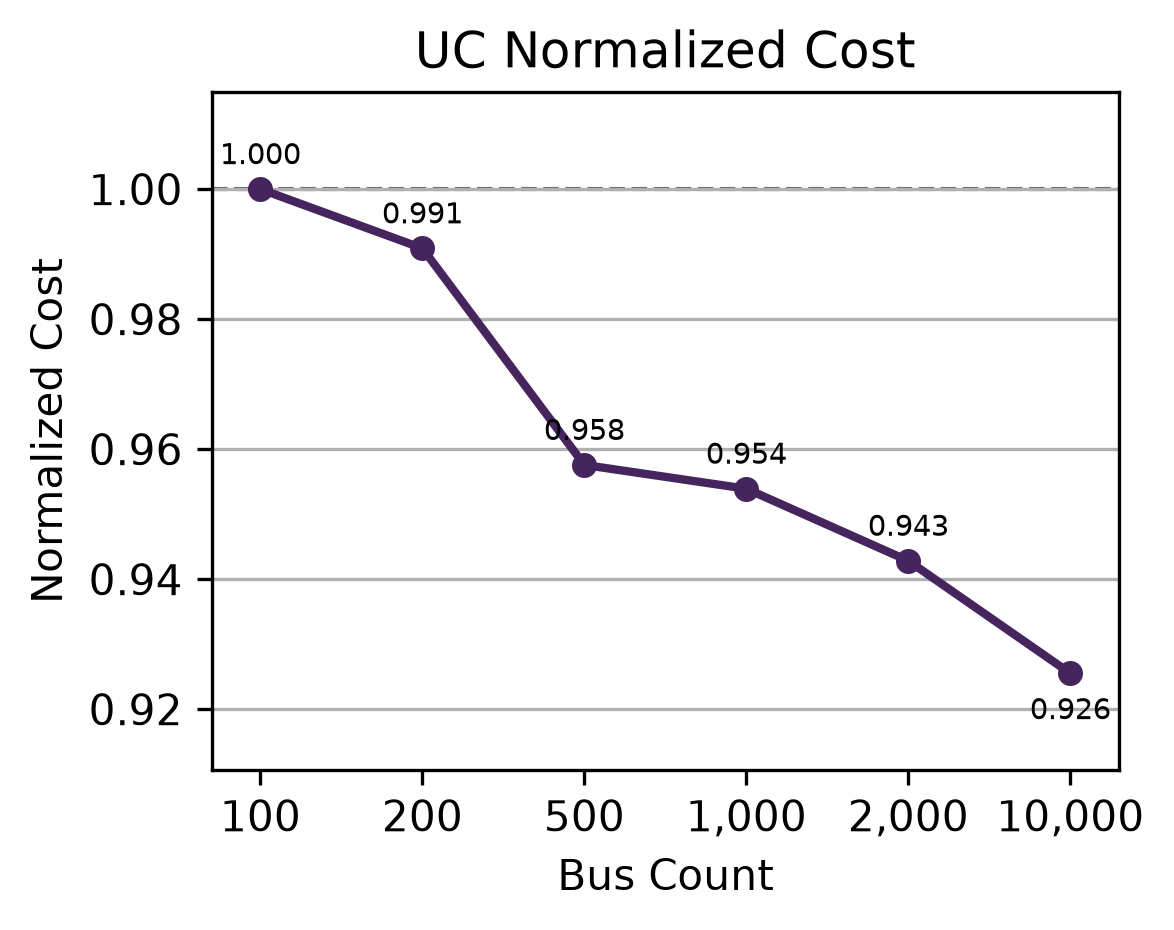}
    \caption{\cref{alg:storage-rruc}'s mean operation cost per bus over a four-day (384 15-min period) horizon, normalized by mean cost per bus for 100 bus case. Optimizing the 10,000 grid bus simultaneously leads to a 7.4\% cost reduction.}
    \label{fig:UC_quasi_cost}
\end{figure}

\cref{fig:UC_quasi_cost} shows the average cost per time period of the UC problem. The normalized cost steadily decreases as the problem grows, showing that optimizing over larger number of buses, generators, and batteries can greatly reduce cost. Optimizing over 10,000 buses reduces costs by about 7.4\% compared to optimizing each 100 bus system separately.

\section{Conclusions}
\label{sec:conclusion}
We extended RRUC for UC without storage~\cite{brody2026fastrruc} to TSD-RRUC for UC with storage, using a TSD of the problem. We used a fast zonal look-ahead plan that fully time couples batteries (but temporarily ignores transmission constraints) to get a target for total zonal battery charge or discharge. We then linked it to a single-period AC UC problem through a marginal price and a net battery charge or discharge target. This keeps all subproblems continuous (even though the original problem has binary variables).

We showed that TSD-ED achieves errors of about 0.1\% from fully time-coupled ED. Combined with our previous results of RRUC's tighter accuracy~\cite{regev2026economicUC,brody2026fastrruc}, this supports a hypothesis that the combined heuristic of time-space decomposition and relaxing and rounding is accurate.

 Our experiments of systems of 100 to 2{,}000 buses show that TSD-ED scales as $\mathcal{O}(n^{1.15})$, compared to fully time coupled ED which scale as $\mathcal{O}(n^{1.66})$. In our tests of 100 to 10{,}000 buses, TSD-RRUC scales as $\mathcal{O}(n^{1.31})$, but we are unaware of a benchmark for comparison. It can solve systems with 10,000 buses, 5,400 generators, and 2,000 batteries in about six minutes on average. This makes it suitable for real-time planning on ISO-size systems.

% Additional ways to analytically or experimentally verify TSD-RRUC's accuracy would make it more usable. Each zone is independent in \cref{eq:target}, so they can be solved in parallel.

\section*{Acknowledgments}
This work was supported in part by the
U.S. Department of Energy, Office of Science, Office of
Workforce Development for Teachers and Scientists (WDTS)
under the Science Undergraduate Laboratory Internships
(SULI) and Community College Internships (CCI) program at Oak Ridge National Laboratory, administered by the Oak Ridge Institute for Science and Education. We thank Slaven Pele\v{s} for writing the proposal that sponsored this work.

\bibliographystyle{IEEEtran}
\bibliography{bibfile}
\end{document}